\documentclass[12pt,a4paper]{article}

\usepackage{setspace}
\usepackage[OT1]{fontenc}
\usepackage{amssymb,amsthm,amsmath,mathabx}

\newcommand\bR{{\mathbb R}}

\newcommand\cov{\mathrm{Cov}}

\newcommand\pto{\stackrel p\rightarrow}

\newcommand\e{{\mathbb E}}
\newcommand\p{{\mathbb P}}
\newcommand\var{{\rm Var}}
\newcommand\tr{{\mathrm{tr}}}

\newtheorem{theorem}{Theorem}[section]%
\newtheorem{lemma}[theorem]{Lemma}%
\newtheorem{proposition}[theorem]{Proposition}%

\begin{document}

\begin{center}
\Large The existence of Riemann-Stieltjes integrals  with applications to  fractional Brownian motion.  
\end{center}
\begin{center}
\large Pavel~Yaskov\footnote{Steklov Mathematical Institute, Russia\\
 e-mail: yaskov@mi.ras.ru\\Supported
    by RNF grant 14-21-00162 from the Russian Science Foundation.}
 \end{center}

\begin{abstract} We derive general sufficient conditions for the existence of Riemann-Stieltjes integrals $\int_a^b Yd X$. Our results extend  the classical conditions of L.C.Young and improve some recent results that deal with  integrals involving a fractional Brownian motion with the Hurst index $H>1/2$.
\end{abstract}

\begin{center}
\end{center}
\section{Introduction}

The paper contributes to the theory of stochastic integration w.r.t. stochastic processes that are not semimartingales. We study conditions under which the integral $\int_a^b YdX$ exists as a limit of forward or ordinary  Riemann-Stieltjes  sums. 

Deep results concerning  Riemann-Stieltjes integrals belong to Young \cite{Y1}, \cite{Y2}. In \cite{Y2}, Young showed  that the convergence of
\begin{equation}\label{phps}
\sum_{n=1}^\infty \phi^{-1}(1/n)\psi^{-1}(1/n)\quad\text{or, equivalentely,}\quad \int_0^1 \frac{\phi^{-1}(u)\psi^{-1}(u)}{u^2}\,du
\end{equation}
implied the existence of the Riemann-Stieltjes integral $\int_a^b fdg$ for any  $f,g$ with no common discontinuities and of bounded  $\phi$- and $\psi$-variations on $[a,b]$, respectively. Here $\phi,\psi\in \mathcal V$ with  $\mathcal V$ being the class  of continuous strictly increasing functions $\varphi$ on $\bR_+$ with $\varphi(\bR_+) =\bR_+$,   and the $\varphi$-variation  of   a function $f=f(t)$ on $[a,b]$ is defined by 
\[V_\varphi(f;[a,b])=\sup \sum_{i=1}^n\varphi(|f(t_i)-f(t_{i-1})|),\]
where the supremum is taken over all partitions $a=t_0<t_1<\ldots<t_n=b$. If $\phi(x)=x^p$ and $\psi(x)=x^q$ for some $p,q>0,$ we get $p$- and $q$-variations, respectively, and \eqref{phps} holds whenever $1/p+1/q>1$. 

Dudley and Norvai\v{s}a \cite{D} provided an extensive list of results concerning $\varphi$-variations of stochastic processes and described different applications of Young's and related results in probability and statistics (see Chapter 12). In \cite{R}, Ruzmaikina studied Riemann-Stieltjes integrals with H\"older continuous functions, a special case covered by Young's results, and used upper bounds on such integrals to prove the existence and uniqueness of solutions of ordinary differential equations with H\"older continuous forcing. Related problems were also studied by Lyons and his collaborators in the context of the rough path theory (e.g., see \cite{L}).

However, the Young theorem   does not cover many cases of interest in stochastic calculus, in particular, certain integrals w.r.t. a fractional Brownian motion  $B_H=B_H(t)$, $t\geqslant 0,$ with the Hurst index $H>1/2$.  The most prominent example is  $\int_0^1 I(B_H(t)>0)dB_H(t)$. This integral arises in non-semimartingale models of the stock market (see \cite{A}). It is hard to check its existence by referring to \eqref{e1}, since, by the self-similarity of $B_H$, the process  $X(t)=I(B_H(t)>0)$ has unbounded $\varphi$-variation on $[0,1]$ for any $\varphi\in \mathcal V$ in any  sense (see Proposition \ref{p2} below).\footnote{Note, however, that the Young theorem guarantees the existence of the Riemann-Stieltjes integral $\int_c^1 I(B_H(t)>0)dB_H(t)$  (e.g., in $L_1$)  for each $c\in(0,1).$} The latter motivated a number of studies \cite{A},  \cite{C}, \cite{S},  where similar stochastic integrals were defined as generalized Lebesgue integrals introduced by Z\"ahle \cite{Z} and further studied by  Nualart and R\u{a}\c{s}canu \cite{N}. It was also shown in \cite{A}, \cite{C},  \cite{S} that these integrals coincide with the Riemann-Stieljes integrals under certain assumptions. 

The purpose of the present paper is to extend condition \eqref{phps} in a simple and natural way to cover the cases of interest involving a fractional Brownian motion or similar stochastic  processes. Our arguments are close to those in ordinary calculus as in Chapter 3 of Dudley and Norvaisa \cite{D}.  In particular, we don't use any fractional derivatives (as in \cite{A}, \cite{C}, \cite{S}).

The paper is structured as follows. Main results are presented in Section 2. Section 3 deals with applications. Proofs are given in Section 4 and an Appendix.

\section{Main results}
Let $X=X(t)$ and $Y=Y(t)$ be stochastic processes on $[a,b]$ and let there exist a function $f=f(s,t),$ $a\leqslant s\leqslant t \leqslant b,$ such that, for all $s<u<t$, 
\begin{equation}
\label{e1}
\|(Y(u)-Y(s))\,(X(t)-X(u))\|_p\leqslant f(s,t),
\end{equation}
 where $p\geqslant 1$ and $\|\xi\|_p$ is the $L_p$-norm. In principal, one can use any other norm or a subadditive functional  on the space of random variables. If  $X=X(t)$ and $Y=Y(t)$ are deterministic, then  $\|\cdot\|_p$ should be replaced by  the absolute value $|\cdot|$. Let also  $f=f(s,t)$  be  non-increasing  in $s$ and non-decreasing in $t$ and $f(t,t)=0$.

Let us introduce the key quantity 
\[I_f(a,b)=\int_{a}^b\int_s^b \frac{f(s,t)}{(t-s)^2}\,dt\,ds+\int_{a}^b\frac{f(a,t)}{t-a}\,dt+
\int_{a}^b\frac{f(s,b)}{b-s}\,ds+f(a,b)\]
for  $f=f(s,t)$ and $a<b$. For a given partition $a=t_0<t_1<\ldots<t_n=b$, define a forward integral sum by
\[\sum_{i=1}^nY(t_{i-1})\Delta X(t_{i})\quad \text{with}\quad \Delta X(t_{i})=X(t_{i})-X(t_{i-1}).\]
Let also  $\max\{t_i-t_{i-1}:1\leqslant i \leqslant n\}$ be the mesh of a partition $\{t_i\}_{i=0}^n$. In what follows all limits of integral sums will be taken over  all sequences of  partitions with mesh tending to zero. 

\begin{theorem} \label{t1}
Let  \eqref{e1} hold and $I_f(a,b)<\infty.$ 
 Then $\int_a^b YdX$ exists as  a limit of  forward integral sums  in $L_p$.
Moreover, for  any   $a=t_0<t_1<\ldots<t_n=b$,
\begin{equation}\label{e2}
\Big\|\int_a^bYdX-\sum_{i=1}^nY(t_{i-1})\Delta X(t_{i})\Big\|_p\leqslant 92\, I_g(t_0,t_n), 
\end{equation}
where $g(s,t)=f(s,t)I(t<s+3d),$ $s\leqslant t,$ and $d$ is the mesh of  $\{t_i\}_{i=0}^n.$ 
\end{theorem}

The proof is based on the following inequality for forward integral sums. 

\begin{theorem}\label{t2}
Let  \eqref{e1} hold. Then, for any partition  $a\leqslant t_0<\ldots<t_{n}\leqslant b,$ 
\begin{equation}\label{ineq}
\Big\|\sum_{i=1}^n(Y(t_{i-1})-Y(t_0))\Delta X(t_{i})\Big\|_p\leqslant 8I_f(t_0,t_n).\end{equation}
\end{theorem}
\noindent{\bf Remark 1.} Instead of forward integrals, one can consider  Riemann-Stieltjes integrals $\int_a^b YdX$    defined as a limit of integral sums $\sum_{i=1}^nY(s_i)\Delta X(t_i)$ with $s_i\in[t_{i-1},t_i].$ All above theorems still hold in this case if $X$ and $Y$ satisfy 
\begin{equation}\label{e3}
 \|(Y(u)-Y(s))\,(X(t)-X(w))\|_p\leqslant f(s,t)
\end{equation} 
 for all $s<t$ and $u,w\in(s,t).$ This can be seen by inspecting the proofs.

\noindent{\bf Remark 2.}  The $L_p$-norm $\|\cdot\|_p,$ $p\geqslant 1,$ in Theorems \ref{t1} and \ref{t2}  can be replaced by $\|\cdot\|_0$ defined by
$\|\xi\|_0=\e[|\xi|/(1+|\xi|)]$ for each random variable $\xi$.  The latter is a subadditive functional determining convergence in probability, i.e.  $\|\xi+\eta\|_0\leqslant \|\xi\|_0+\|\eta\|_0$ for any $\xi,\eta$ and
\[\xi_n\pto\xi\text{ for some $\xi$ iff }\|\xi_n-\xi_m\|_0\to0,\;m,n\to\infty.\]
These are the only  properties of $\|\cdot\|_p$ used in the proofs of Theorems \ref{t1} and \ref{t2}. Note also that, by Jensen's inequality, for all $p\geqslant 1$ and $\xi,$
\[\|\xi\|_0\leqslant \frac{\e|\xi|}{1+\e|\xi|}\leqslant \frac{\|\xi\|_p}{1+\|\xi\|_p}\]

The "norm" $\|\cdot\|_0$ can be  more useful than $\|\cdot\|_p$. In the typical case with $f(s,t)=C(t-s)^\alpha s^{-\beta}$ in  \eqref{e1}   for some $C>0,$ $p\geqslant 1,$ $\alpha>1,$ and $\beta\in(1,\alpha)$, we have  $I_f(0,1)=\infty$. But it can be checked that $I_g(0,1)<\infty$ for
\[g(s,t)=\frac{f(s,t)}{1+f(s,t)}=\frac{Ct^{\alpha-\beta}(1-s/t)^\alpha}{(s/t)^\beta +Ct^{\alpha-\beta}(1-s/t)^{\alpha} }.\]

\section{Applications}

Let us first discuss the relation of our results to the  Young theorem stated in the Introduction. To do it, we need some preliminary facts.  

Note  that any forward (or Riemann-Stieltjes) integral sum over $[a,b]$ stays the same when changing the time $t\to s=\varphi(t)$ for an increasing continuous function $\varphi$. That is,
\[\sum_{i=1}^nY(t_{i-1})\Delta X(t_{i})=\sum_{i=1}^nY(\varphi^{-1}(s_{i-1}))\Delta X(\varphi^{-1}(s_{i}))\]
whenever $s_i=\varphi(t_i),$ $i=0,\ldots,n.$ Thus, Theorem  \ref{t1} will still hold  if we change the time $s=\varphi(t)$ and replace  $f=f(s,t)$  by $f_\varphi=f(\varphi^{-1}(s),\varphi^{-1}(t))$ and $[a,b]$ by  $[\varphi(a),\varphi(b)]$. 

Fix $[a,b]\subset \bR.$ As above, let  $\mathcal V$ be the class of continuous increasing functions $\varphi$ on $\bR_+$  with $\varphi(\bR_+)=\bR_+$ and, for $p\geqslant 1$ and $\varphi\in\mathcal V,$ define \[\|X\|_{\varphi,p}=\inf\{C>0:V_{\varphi,p}(X/C;[a,b])\leqslant 1\},\]  where $V_{\varphi,p}$ is defined as $V_\varphi$ in the Introduction  with $|\cdot|$ replaced by $\|\cdot\|_p.$

\begin{proposition}\label{p1}
For $q,r>0$, let $X=X(t)$ and $Y=Y(t)$ be stochastic  processes continuous on $[a,b]$ in $L_q$ and $L_r$, respectively. If $p=(1/q+1/r)^{-1}$ and $\|X\|_{\phi,q},\|Y\|_{\psi,r}<\infty,$ then there is an increasing continuous function $\varphi$ on $[a,b]$ such that \eqref{e3} holds for \[f(s,t)= \|X\|_{\phi,q} \|Y\|_{\psi,r}\,\phi^{-1}(\varphi(t)-\varphi(s))\psi^{-1}(\varphi(t)-\varphi(s)). \] Moreover, for $g=f_{\varphi},$
\[I_{g}(\varphi(a),\varphi(b))\leqslant  \|X\|_{\phi,q} \|Y\|_{\psi,r}\Big(\int_0^\infty \frac{\phi^{-1}(u)\psi^{-1}(u)}{u^2}du+\phi^{-1}(3)\psi^{-1}(3)\Big).\]
\end{proposition}
Now, it can be readily seen that Proposition \ref{p1} and Theorem \ref{t1} extend the Young theorem (see the Introduction) at least in the case of continuous processes. To ease the presentation, we won't consider the case of processes with no common discontinuities. 

Let us apply our results to integrals w.r.t. 
a fractional Brownian motion   $B_H=B_H(t),$ $t\geqslant 0,$ with the Hurst index $H\in(1/2,1)$. Namely, we are interested in the existence of the Riemann-Stieltjes integrals $\int_0^TF(B_H)dB_H$ for  $T>0$ and a given function $F=F(x).$ 

In general, the Young theorem doesn't guarantee the existence of such integrals. In the special case $F(x)=I(x>0),$ the latter follows from the next proposition.
\begin{proposition}\label{p2}
Let $\varphi\in \mathcal V,$  $H\in(1/2,1)$, and $X(t)=I(B_H(t)>0)$, $t\geqslant 0.$ Then $X=X(t)$  has infinite $\varphi$-variation on $[0,1]$ a.s. and in any norm $\|\cdot\|$ on the space of random variables $($$\|\xi\|$  depends only on the distribution of $\xi$$)$.
\end{proposition}

We will further consider only  non-decreasing $F$.   For any such $F$, there is a $\sigma$-finite measure $\mu$ on the Borel $\sigma$-algebra in $\bR$ such that \[F(x)=F(0)+\int_0^x \mu(du)\text{\quad at any continuity point $x\in\bR$ of $F$.}\] We also have the following.

\begin{proposition}\label{p}  If $F=F(x)$ is non-decreasing on $\bR$, then the following conditions are equivalent:
$(a)$  $\int_{-\infty}^\infty \exp\{-\varepsilon u^2\}\mu(du)<\infty$ for all $\varepsilon>0$,  $(b)$  $\e |F(B_H(t))|<\infty$ for all $t>0$.
\end{proposition}

\begin{theorem}\label{t3} Let $(a)$ in Proposition \ref{p} hold, $T>0$, $H\in(1/2,1)$, and $p\in[1,H/(1-H))$. Then the Riemann-Stieltjes integral $I=\int_0^TF(B_H)dB_H$ exists as a limit in $\|\cdot\|_{p}$ and 
\[\int_0^{B_H(T)}F(x)\,dx=\int_0^T F(B_H)dB_H\text{ a.s.}\] Moreover, if  $\{t_{i}\}_{i=0}^n$ is a partition of $[0,T]$ with the mesh $d>0,$ then
\[\Big\|I-\sum_{i=1}^nF(B_H(t_{i-1}))\Delta B_H(t_i)\Big\|_p\leqslant CT^{1-H/q}d^{H/q+H-1}\]
for some $C=C(p,q,T,H,F)$ and any $q\in(p,H/(1-H))$. 
\end{theorem}
 Taking $p=1$ and $q$ close enough to $p$ in Theorem \ref{t3},  we derive 
\[\sum_{i=1}^nF(B_H(t_{i-1}))\Delta B_H(t_i)=\int_0^TF(B_H)dB_H +O_p(d^{2H-1-\varepsilon}),\quad d\to 0,\]
when $T,\varepsilon>0 $ and $H>1/2$ are fixed, and $d$ is the mesh of $\{t_i\}_{i=0}^n$.

Theorem \ref{t3} improves Theorem 3.1 in \cite{AV}, where the uniform partitions are considered and the established rate of convergence is $O_p(d^{H-1/2-\varepsilon})$. Related results are Theorem 4.1 in \cite{A}, Theorem 3.3 in \cite{S}, and Theorem 4.1 in \cite{C}. These theorems deal with similar integrals defined in a pathwise sense. However, as it is discussed in \cite{C}, the proofs in  \cite{A} and \cite{S}  contain  some gaps that are covered in \cite{C}. Another related result  is Theorem 2.2 in \cite{SZ}, where discrete approximations (in $L_p$) of forward integrals $\int_0^t XdB_H$ are studied for $X=X(t)$ of finite integral $q$-variation (for details, see \cite{SZ}).

In the proof of Theorem \ref{t3}, the following inequality plays a key role. 
\begin{proposition}\label{ll1}  There exists an absolute constant $C>0$ such that \[\p(B_H(s)<v<B_H(t))\leqslant C\exp\{-v^2/(4t^{2H})\}\frac{(t-s)^H}{t^H} \]
for all  $H\in[1/2,1)$ and  $0\leqslant s<t$.
\end{proposition}

Proposition \ref{ll1} improves Lemma 4 in \cite{CN} and Lemma 3.2 in \cite{MSV}  (see also Lemma 3.1 in \cite{AV}).  Under the conditions of Theorem \ref{t3}, it allows to show that   \eqref{e1} holds for    $f(s,t)=C(t-s)^{H/q+H}t^{-H/q},$ $s\leqslant t,$
$X(t)=B_H(t)$ and $Y(t)=F(B_H(t)),$ $t\geqslant 0,$ where $C>0$ do not depend on $(s,t)$. In fact, the latter is true for any $X=X(t)$ such that $\|X(t)-X(s)\|_p\leqslant C_p|t-s|^H$ for all $s,t\in[0,T]$ and $p\geqslant 1$ (e.g., $X(t)=e^{B_H(t)}$, see Lemma 3.6 in \cite{T}).

\section{Proofs}
\noindent{\bf Proof of Theorem \ref{t1}.} Let $t_0<\ldots< t_n$ and $s_0<\ldots<s_m$ be partitions of $[a,b]$.  Let $q_0<\ldots<q_l$ be a partition of $[a,b]$ obtained by taking the union of $\{s_i\}_{i=0}^m$ and $\{t_i\}_{i=0}^n$. Here $l\leqslant n+m$. Then $I= \sum_{j=1}^{l}Y(q_{j-1})\Delta X(q_j)$ can be written as
\[\sum_{i=1}^{n}Y(t_{i-1})\Delta X(t_i) + \sum_{i=1}^{n}\sum_{j:\,t_{i-1}\leqslant q_{j-1}<t_i }(Y(q_{j-1})-Y(t_{i-1}))\Delta X(q_j) .\]
Hence, by Theorem \ref{t2},
\begin{equation}\label{DIc1-0}\Bigl\|I-\sum_{i=1}^{n}Y(t_{i-1})\Delta X(t_i)\Bigr\|_p\leqslant
8\sum_{i=1}^n I_f(t_{i-1},t_i).\end{equation} Replacing $\{t_i\}_{i=0}^n$ by  $\{s_i\}_{i=0}^m$,
we get a  similar bound with $I_f(s_{i-1},s_i)$. Hence, the integral sums over $\{t_i\}_{i=0}^n$ and  $\{s_i\}_{i=0}^m$ differ no more than 
\begin{equation*}
8\sum_{i=1}^nI_f(t_{i-1},t_i)+8\sum_{i=1}^mI_f(s_{i-1},s_{i})
\end{equation*}
in the $L_p$-norm. Thus, $\int_a^bY dX $ exists if \begin{equation}
\label{dd}\text{$\sum_{i=1}^nI_f(t_{i-1},t_i)\to 0$ as the mesh of a partition $\{t_i\}_{i=0}^n$ tends to 0.}
\end{equation} To prove it, we  need the following lemma (for the proof, see the Appendix).
\begin{lemma}\label{l1}
Let $a\leqslant t_0<\ldots<t_n\leqslant b$ and $d=\max\{t_i-t_{i-1}:1\leqslant i\leqslant n\}$. If $d<1/6$ and $g(s,t)=f(s,t)I( t < s+ 3d)$, $s\leqslant t,$  then 
\[\sum_{i=1}^nI_f(t_{i-1},t_i)\leqslant
92\,I_g(t_0,t_n).\] 
\end{lemma}
If $I_f(a,b)$ is finite, then \eqref{dd} follows from Lemma \ref{l1} and the absolute continuity of Lebesgue's integral. We obtain the bound \eqref{e2} by taking a limit over $\{s_i\}_{i=0}^m$ in 
\eqref{DIc1-0} and applying Lemma \ref{l1}. This proves the theorem. $\square$

\noindent{\bf Proof of Theorem \ref{t2}.}   For simplicity, let $[a,b]\subseteq [0,1]$ and $T=\{t_i\}_{i=0}^n .$ Define $j$ by $1/2\in [t_{j},t_{j+1})$ if such $j$ exists otherwise put $j=n$ for $t_n<1/2$ and $j=0$ for $t_0>1/2$. Denoting the integral sum in \eqref{ineq} by $I$, we can write
\[I= \sum_{i=1}^{j}(Y(t_{i-1})-Y(t_0))\Delta X(t_{i}) +\sum_{i=j+2}^{n}(Y(t_{i-1})-Y(t_{j+1}))\Delta X(t_{i})+R,\]
where any sum  over the empty set is zero  and the remainder is defined by
\[R=[Y(t_{j})-Y(t_{0})]\,[X(t_{j+1})- X(t_{j})]+[Y(t_{j+1})-Y(t_0)]\,[X(t_{n})- X(t_{j+1})].
\]
By \eqref{e1} and the monotonicity of $f$, $\|R\|_p\leqslant 2f(0,1)$. By construction, $T_0=\{t_0,\ldots,t_{j}\}\subseteq[0,1/2]$ and $T_{1}=\{t_{j+1},\ldots,t_{n}\}\subseteq[1/2,1]$, where $\{t_i,\ldots,t_k\}$ is empty if $i\geqslant k$. By the same arguments, we can write 
the integral sum  over each $T_i$ as a sum of two integral sums over some $T_{i
,0}\subseteq[(2i)/4, (2i+1)/4]$ and $T_{i,1}\in [(2i+1)/4, (2i+2)/4]$ and a reminder $R_i$ such that
$\|R_1\|_p+\|R_2\|_p\leqslant 2f(0,1/2)+2f(1/2,1)$. 

 Continuing this procedure until  each integral sum is over a two-point set or over the empty set, we get, for $s_i^{n}=i/2^{n}$,
\[\|I\|_p\leqslant 2f(0,1)+2\sum_{n\geqslant 1}
\sum_{i=1}^{2^n}f\bigl(s_{i-1}^{n},s_{i}^{n}\bigr),\]
since such integral sums are zeros.  By the monotonicity of $f$, 
\begin{align*}
\sum_{n\geqslant1}
\bigl(f\bigl(0,s_{1}^{n}\bigr)+f\bigl(s_{2^n-1}^{n},1\bigr)\bigr)&\leqslant \int_0^\infty [f(0,2^{-u})+
f(1-2^{-u},1)]\,du\\&\leqslant\frac{1}{\ln 2}\int_0^1\bigg(\frac{f(0,x)}{x}\,dx+
\frac{f(x,1)}{1-x}\bigg)\,dx,
\end{align*}
\begin{align*}
\sum_{n\geqslant2}\sum_{i=2}^{2^n-1}f\bigl(s_{i-1}^{n},s_{i}^{n}\bigr)&=\sum_{n\geqslant2}
\sum_{i=2}^{2^n-1}\frac{1}{\ln(4/3)}\iint\limits_{A_i^n}\frac{
f\bigl(s_{i-1}^{n},s_{i}^{n}\bigr)}{(s-t)^2}\,dt\,ds\\&\leqslant
4\int_0^1\int_s^1\frac{f(s,t)}{(s-t)^2}\,dt\,ds,
\end{align*}
for $A_i^n=\bigl[s_{i-2}^n,s_{i-1}^n\bigr]\times\bigl[s_i^n,s_{i+1}^n\bigr).$ Here the last inequality holds since $A_i^n$   does not intersect over $k$ for each  $n$. Combining the above bounds we get the desired inequality. Q.e.d.

\noindent{\bf Proof of Proposition \ref{p1}.} Assume w.l.o.g. that $\|X\|_{\phi,q}=1$ and  $\|Y\|_{\psi,r}=1$, otherwise replace $(X,Y)$ by  $(X/\|X\|_{\phi,q},Y/\|Y\|_{\psi,r})$. By definition, $V_{\phi,q}$ (as well as $V_{\psi,r}$) is superadditive, i.e., for all $a\leqslant s<t\leqslant b$,
\[V_{\phi,q}(X;[a,t])\geqslant V_{\phi,q}(X;[a,s])+V_{\phi,q}(X;[s,t]).\]
Hence, by H\"older inequality, for all $s<t$ and  $u,w\in (s,t),$
\begin{align}\label{e4}
\|(Y(u)-Y(s))\,&(X(t)-X(w))\|_p\leqslant \|Y(u)-Y(s)\|_q\,\|X(t)-X(w)\|_r\nonumber\\
&\leqslant  \psi^{-1}(\psi(\|Y(u)-Y(s)\|_q))\,\phi^{-1}(\phi(\|X(t)-X(w)\|_r))\nonumber\\
&\leqslant  \psi^{-1}(V_{\psi,r}( Y;[s,t]))\,\phi^{-1}(V_{\phi,q}( X;[s,t]) )\nonumber\\
&\leqslant  \psi^{-1}(\varphi(t)-\varphi(s))\,\phi^{-1}(\varphi(t)-\varphi(s))=:f(s,t),
\end{align}
where, for any  $a\leqslant t\leqslant  b,$   $\varphi(t)=\varphi_1(t)+\varphi_2(t)+(t-a)/(b-a)$,
\[\text{$\varphi_1(t)=V_{\psi,r}(Y;[a,t])$  and $\varphi_2(t)=V_{\phi,q}(X;[a,t])$.}\] Since $X=X(t)$ and $Y=Y(t)$ are continuous in $L_q$ and $L_r$ respectively, $\varphi$ is a continuous strictly increasing function.  Moreover, $\varphi(a)=0$ and $\varphi(b)\leqslant 3$ (since $\|X\|_{\phi,q}=1$ and  $\|Y\|_{\psi,r}=1$). Set  \[g(s,t)=f(\varphi^{-1}(s),\varphi^{-1}(t))=\phi^{-1}(t-s)\psi^{-1}(t-s),\] we also have that
\[ \int_0^3\int_s^3 \frac{\phi^{-1}(t-s)\psi^{-1}(t-s)}{(t-s)^2}\,dt\,ds\leqslant 3\int_0^3\frac{\phi^{-1}(u)\psi^{-1}(u)}{u^2}\,du=:I,\]
\[ \int_0^3\frac{\phi^{-1}(t)\psi^{-1}(t)}{t}\,dt\leqslant I,\quad  \int_0^3\frac{\phi^{-1}(3-s)\psi^{-1}(3-s)}{3-s}\,dt\leqslant I.\]
We finish the proof by noting that $I_g(\varphi(a),\varphi(b))\leqslant 9I+\phi^{-1}(3)\psi^{-1}(3)
.$ Q.e.d.

\noindent{\bf Proof of Proposition \ref{p2}.}  By the self-similarity of $B_H,$
\[\Delta_n=I(B_H(2^{-n})>0)-I(B_H(2^{-n-1})>0),\quad n=0,1,2,\ldots,\]  are identically distributed and non-degenerate. For any fixed $\varphi\in\mathcal V,$ $\varphi(x)>0$ when $x>0$ and $\varphi(\infty)=\infty.$ Therefore, if 
the norm $\|\xi\|$ depends only on the distribution of a random variable $\xi$, then
\[\sum_{n=0}^\infty\varphi(\|\Delta_n\|)=\infty.\]
Thus,  $X(t)=I(B_H(t)>0)$ has unbounded $\varphi$-variation on $[0,1] $ in $\|\cdot\|.$

Let us also  show that $\sum_{n=0}^\infty\varphi(|\Delta_n|)=\infty$ a.s., 
i.e. $X=X(t)$ has  infinite $\varphi$-variation on $[0,1]$ a.s. 
To see it, note that, for $H>1/2$,
\begin{align} \cov(B_H(t),B_H(s))=&\frac{t^{2H}-(t-s)^{2H}+s^{2H}}{ 2}
\leqslant  \label{a}Ht^{2H-1} s +\frac{s^{2H}}{ 2 }=o(s^H)
\end{align}
uniformly in $t\in[t_0,1]$ for $t_0>0$ as $s\to 0.$ In addition, by the self-similarity of $B_H,$ the spectral norm 
$\|\var(\xi_t)^{-1/2}\|=Ct^{-H}$ for all $t>0$, $\xi_t=(B_t,B_{t/2}),$ and $C=\|\var(\xi_1)^{-1/2}\|.$ Hence, by \eqref{a},  
\[\rho(s,t)=\|\var(\xi_t)^{-1/2}\cov(\xi_t,\xi_s )\var(\xi_s)^{-1/2}\|\leqslant C^2s^{-H}t^{-H}\|\cov(\xi_t,\xi_s)\|\to 0\] uniformly in $t\in[t_0,1]$ for $t_0>0$ as $s\to 0,$
where we have used the fact that \[\|A\|\leqslant \sqrt{\tr(A^\top A)}\leqslant 4\max_{1\leqslant i,j\leqslant 2}|a_{ij}|\] for any $2\times 2$ matrix $A=(a_{ij})_{i,j=1}^2.$ By the Gebelein inequality for Gaussian vectors in $\bR^2$ (see Theorem 3.4 and formula (3.1)  in \cite{V}),
\[\mathrm{Corr}(I(|\Delta_{n }|=1),I(|\Delta_{m}|=1))\leqslant C\rho(2^{-n},2^{-m}),\quad m,n\geqslant 1,\]
for an absolute constant $C>0$. Hence, there exists a sequence $(n_k)_{k=1}^\infty $  such that $n_k\uparrow \infty$ as $k\uparrow \infty$ and each $n_{l},$ $l>1,$ is chosen in a way  that
\[\p(|\Delta_{n_k}|=1,|\Delta_{n_{l}}|=1)\leqslant  \p(|\Delta_{n_k}|=1)\p(|\Delta_{n_{l }}|=1)+\frac{1}{2^{k+l}},\quad k<l,\]
given $n_1,\ldots,n_{l-1}.$ By the Erd\"os-Renyi theorem, if $\sum_{k=1}^\infty \p(A_k)=\infty,$ then 
\[ \p( A_k\;\text{i.o.})\geqslant\varlimsup_{n\to\infty}\frac{(\sum_{i=1}^n\p(A_i))^2}{\sum_{i,j=1}^n\p(A_iA_j)},\] where i.o. means infinitely often. Taking $A_k=\{|\Delta_{n_k}|=1\}$, we get 
\[ \p( \Delta_{n_k}=1\;\text{i.o.})\geqslant\varlimsup_{k\to\infty}\frac{(np)^2}{np+n^2p^2+\sum_{i,j=1}^\infty 2^{-i-j}}=1 \] for $p=\p(|\Delta_1|=1)>0.$  This yields $\sum_{n=1}^\infty\varphi(|\Delta_n|)=\infty$ a.s. Q.e.d.

\noindent{\bf Proof of Proposition \ref{p}.} Assume w.l.o.g. that $F(0)=0$. As a result, $F(x)=\int_0^x\mu(du)$ for any continuity point $x$ of $F.$ Write further $B_t$ instead of $B_H(t)$. Fix $t>0.$ We have
\[\e |F(B_t)|=\e F(B_t) I(B_t>0)-\e F(B_t) I(B_t<0).\] By the Fubini-Tonelli theorem,
\begin{align*}
\e F(B_t) I(B_t>0)&=\e \int_0^\infty I(B_t\geqslant u)\mu(du) =\int_0^\infty\p(\xi>ut^{-H})\mu(du),
\end{align*}
where  $\xi\sim\mathcal N(0,1)$. It is straightforward to check (see \eqref{Phi}) that  \[C\exp\{-u^2t^{-2H}\}\leqslant \p(\xi>ut^{-H})\leqslant \exp\{-u^2t^{-2H}/2\} \] for all $u>0$ and some $C>0$. The term $\e F(B_t) I(B_t<0)$ can be analysed similarly. Taking $\varepsilon=t^{-2H}$ or $\varepsilon=t^{-2H}/2$, we see that $(a)$ and $(b)$ are equivalent. Q.e.d.

\noindent{\bf Proof of Theorem \ref{t3}.} Write  $B_t$ instead of $B_H(t)$ and define $Y(t)=F(B_t)$ and $X(t)=B_t$, $t\geqslant 0.$ Fix $s,u,w,t$ with $s<t$ and $u,w\in(s,t)$. 
For any $r\geqslant 1$, there is $C_r>0$ such that $\|X(t)-X(w)\|_r= C_r(t-w)^{H}$.
For fixed $q>p\geqslant 1,$ define $r>p$ from $1/r+1/q=1/p$. By H\"older's inequality, 
\begin{align*}
\|(Y(u)-Y(s))\,(X(t)-X(w))\|_p&\leqslant \|Y(u)-Y(s)\|_q \|X(t)-X(w)\|_r\\
&\leqslant C_r\|Y(u)-Y(s)\|_q (t-s)^H.
\end{align*}
Since $F=F(x)$ is continuous in $x=B_t$ a.s., we have
\[|Y(u)-Y(s)|=\Big|\int_{B_s}^{B_u} \mu(dv)\Big|=\int_\bR [I(B_s<v<B_u)+I(B_s>v>B_u)]\mu(dv).\]
Therefore, $\|Y(u)-Y(s)\|_q\leqslant\|I_1\|_q+\|I_2\|_q$, where
\[I_1=\int_\bR I(B_s<v<B_u)\mu(dv)\quad\text{and}\quad I_2=\int_\bR I(B_u<v<B_s)\mu(dv).\]
Since $(-B_t)_{t\geqslant 0}$ has the same distribution as $(B_t)_{t\geqslant 0}$, we will estimate only $\|I_1\|_q.$ Fix $\varepsilon>0$ and set $A_\varepsilon=\int_\bR \exp\{-\varepsilon v^2\}\mu(dv)$. By  Lyapunov's inequality, 
\begin{align*}
\|I_1\|_q^q&=\e\bigg(\int_\bR I(B_s<v<B_u)\mu(dv)\bigg)^q\\&=A_\varepsilon^q\e\bigg(\int_\bR I(B_s<v<B_u)\exp\{\varepsilon v^2\}\frac{e^{-\varepsilon v^2}\mu(dv)}{A_\varepsilon}\bigg)^q\\&\leqslant A_\varepsilon^q\e\int_\bR I(B_s<v<B_u)\exp\{q\varepsilon v^2\}\frac{e^{-\varepsilon v^2}\mu(dv)}{A_\varepsilon}\\
&\quad=A_\varepsilon^{q-1}\int_\bR \p(B_s<v<B_u)\exp\{(q-1)\varepsilon v^2\}\mu(dv).
\end{align*}
By Proposition \ref{p} and the monotonicity of $g(u)=(u-s)/u$ on $[s,t]$, 
\begin{align*}
\|I_1\|_q^q&\leqslant CA_\varepsilon^{q-1}
\int_\bR \exp\{-v^2/(4u^{2H})+(q-1)\varepsilon v^2\}\mu(dv) \frac{(u-s)^H}{u^{H}}\\
&\leqslant CA_c^{q } \frac{(t-s)^H}{t^{H}},
\end{align*}
where $c=\min\{\varepsilon,\delta\}$ for $\delta=1/(4T^{2H})-(q-1)\varepsilon$ and we take small enough $\varepsilon$ to get $\delta>0.$
Finally, we see that \eqref{e1} holds for \[f(s,t)=K (t-s)^{H/q+H}t^{-H/q}=K(1-s/t)^{H/q} (t-s)^{H} ,\]
where $K$ depends on $T,$ $H,$ $p,$ $q$, and $F$.  The function $f=f(s,t),$ $s\leqslant t,$ increases in $t$ and decreases in $s.$ It easy to see that all integrals in the definition of $I_f(0,T)$
converge when $H/q<1$ and $H/q+H>1.$ 

By Theorem \ref{t1}, $I=\int_0^TF(B_H)dB_H$ exists as a limit in $L_p$. 
Moreover, if $\{t_i\}_{i=0}^n$ is a partition of $[0,T]$ with the mesh $d>0$, then 
\[\Big\|I-\sum_{i=1}^nF(B_H(t_{i-1}))\Delta B_H(t_i)\Big\|_p\leqslant 92 \,I_g(0,T),\]
where $g(s,t)=f(s,t)I(t<s+3d)$, $s\leqslant t$. We have
\begin{align*}
\int\limits_0^T\int\limits_{s}^{s+3d}\frac{(t-s)^{H/q+H} }{(t-s)^2}\frac{dtds}{t^{H/q}}&\leqslant
\int\limits_0^T\frac{dt}{t^{H/q}}\int\limits_0^{3d}u^{H/q+H-2}du=\frac{T^{1-H/q}(3d)^{H/q+H-1}}{(1-H/q)(H/q+H-1)},\\
\int\limits_0^{3d}\frac{t^{H/q+H}}{t}\frac{dt}{t^{H/q}}&\leqslant \int\limits_0^{3d}t^{H -1}dt=\frac{(3d)^{H }}{H }\leqslant \frac{(3d)^{H } (T/d)^{1-H/q}}{H },\\
\int\limits_{T-3d}^{T}\frac{(T-s)^{H/q+H}}{T-s}\frac{ds}{T^{H/q}}&\leqslant \int\limits_{T-3d}^{T}(T-s)^{H-1}ds=\frac{(3d)^{H}}{H }\leqslant\frac{(3d)^{H } (T/d)^{1-H/q}}{H }.\end{align*}
 As a result, 
\[\Big\|I-\sum_{i=1}^nF(B_H(t_{i-1}))\Delta B_H(t_i)\Big\|_p\leqslant 92\,I_g(0,T)\leqslant C T^{1-H/q}d^{H/q+H-1} \]
for  $C=C_0A_c$ and $C_0=C_0(p,q,H)>0.$

Let us show that
\begin{equation}\label{ito}
\int_0^{B_H(T)}F(x)\,dx=\int_0^TF(B_H )\,dB_H \text{ a.s.} 
\end{equation}
First, let $F$ be smooth.  It is well known that  $B_H=B_H(t)$  has zero quadratic variation on $[0,T]$ a.s. when $H>1/2.$ As it was shown in \cite{F},  $\int_0^TF(B_H )\,dB_H $ exists as a limit of forward integral sums a.s. and, as a result, it coincides with the limit in $L_p$. Thus, \eqref{ito} holds for smooth $F$. 

The general case with can be obtained by taking a limit over smooth $F_n$ tending to $F.$ Let $(F_n)_{n=1}^\infty$ be a sequence of smooth non-decreasing functions with $F_n(x)\to F(x)$ in all continuity points $x$ of $F$, then, by the Lebesgue dominated convergence theorem, 
\begin{equation}\label{ci}
\int_0^{B_H(T)}F_n(x)\,dx\to \int_0^{B_H(T)}F(x)\,dx\text{ a.s.}
\end{equation}
To prove \eqref{ito}, we only need to show that 
\begin{equation}
\label{In}I_n=\int_0^TF_n(B_H )\,dB_H\pto \int_0^TF(B_H )\,dB_H=I.
\end{equation}
For any fixed partition  $\{t_i\}_{i=0}^m$ of $[0,T]$ with the mesh $d>0$, \[\sum_{i=2}^mF_n(B_H(t_{i-1}))\Delta B_H(t_i)\pto \sum_{i=2}^mF(B_H(t_{i-1}))\Delta B_H(t_i),\quad n\to\infty,\]
since $F=F(x)$ is continuous in $x\in\{B_H(t_{i})\}_{i=0}^n$ a.s. 

As it is shown above,
\[\Big\|I_n-\sum_{i=2}^mF_n(B_H(t_{i-1}))\Delta B_H(t_i)\Big\|_p\leqslant C_1A_{c}^n  d^{H/q+H-1}+K_p|F_n(0)|d^{H},\] 
where $A_{c}^n=\int_\bR \exp\{-c u^2\}dF_n (u)$ with $c$ given above, $C_1=C_1(p,q,T,H)$, and the last term appears as an estimate of $\|F_n(0)\Delta B_H(t_1)\|_p$.  The same bound holds for $(I_n,F_n,A_c^{n}) $ replaced by $(I,F,A_c).$ These bounds show  that the above integrals can be uniformly approximated by integral sums whenever  $A_c^{n}=O(1)$. By \eqref{ci} and \eqref{In}, the latter will yield   \eqref{ito}.

To finish the proof, we need to choose smooth $F_n$ such that   $A_c^n=O(1)$ and $F_n(x)\to F(x)$ for all continuity points $x$ of $F$.  This is true for
 \[F_n(x)= \int_{0}^{1}\varphi(z) F(x+z/n)\,dz\]
where $\varphi=\varphi(z)$ be a $C^\infty$-density whose support is $[0,1]$. Q.e.d.

\noindent{\bf Proof of Proposition \ref{ll1}.} Fix  $H\in[1/2,1)$ and  $t>s>0$ (the case $s=0$ can be obtained by continuity).  In what follows, write  $B_s$ and $B_t$ instead of $B_H(s)$ and $B_H(t)$. 
Set $\alpha=\cov(B_t,B_s)/\var(B_s)$. Then $Z=B_t-\alpha B_s$ is independent of $B_s$. Moreover, $\var(Z)=\var(B_t)-\var(\alpha B_s)=t^{2H}-\alpha^2 s^{2H}$ and
\begin{equation}\label{ts}
t^{2H}-\alpha^2 s^{2H}=\min_{a\in \bR}\e|B_t-aB_s|^2\leqslant \e|B_t -B_s|^2=(t-s)^{2H}.
\end{equation}

Next, by the convexity of $f(x)=x^{2H}$  on $\bR_+$ (for $H>1/2$),  \[\alpha-1=[(t^{2H}-(t-s)^{2H})-(s^{2H}-0^{2H})]/(2s^{2H})\geqslant 0.\] Since $f(x)=x^{H}$ is concave on $\bR_+$ and $g(x)=(a-x)/(b-x)$ is decreasing over $x$ when $x\in[0,a]$ and $0<a\leqslant b$,
\begin{align}\label{1/alp}
1-\frac 1\alpha=&\frac{t^{2H}-s^{2H}-(t-s)^{2H}}{t^{2H}+s^{2H}-(t-s)^{2H}}\leqslant \frac{t^{2H}-s^{2H}}{t^{2H}+s^{2H}}\leqslant \nonumber\\
&\leqslant \frac{(t^{H}-s^{H})(t^{H}+s^{H})}{(t^{H}+s^{H})^2/2}\leqslant \frac{2(t^{H}-s^{H})}{t^H}\leqslant \frac{2(t-s)^H}{t^H}.
\end{align} 
We will also need the following useful facts. If $\xi\sim \mathcal N(0,1)$ and $\Phi=\Phi(x)$ is its c.d.f., then, for all $x>0,$ $\Phi(-x)=1-\Phi(x)$ and
\begin{equation}\label{Phi}
1-\Phi(x)=\p(\xi>x)\leqslant \e \exp\{x\xi-x^2\}=\exp\{-x^2/2\}.
\end{equation}

First, let $t>4s/3.$ Then  $4^{-1}\leqslant 4^{-H}\leqslant (1-s/t)^H$. If $v>0$, then
\[I=\p(B_s<v<B_t)\leqslant \p(B_t>v)\leqslant\]\[\leqslant  \exp\{-v^2/(2t^{2H})\}\leqslant 4\exp\{- v^2/(2t^{2H})\} \frac{(t-s)^H}{t^H}.\]
If $v<0$, then  
\[I\leqslant \p(B_s<v)\leqslant  \exp\{-v^2/(2s^{2H})\}\leqslant 4 \exp\{-v^2/(2t^{2H})\} \frac{(t-s)^H}{t^H}.\]

Now, let $s<t\leqslant 4s/3$. Then 
\begin{equation}\label{alp}
\alpha=\frac{t^{2H}+s^{2H}-(t-s)^{2H}}{2s^{2H}}\leqslant \frac{(4/3)^{2H}+1}{2}\leqslant \frac{16/9+1}{2}=\frac{25}{18}<\sqrt{2}.
\end{equation} 
By $\alpha>1,$
\begin{align*}
|\p(B_s<v<B_t)&-\p(\alpha B_s<v< B_t)|\leqslant \\
&\leqslant \e|I(B_s<v)-I(B_s<v/\alpha)|=I_1,
\end{align*}
where $I_1=|\p(B_s<v)-\p(B_s<v/\alpha)|.$ Now, if $\varphi =\Phi' $, then there is $C>0$ such that  $|x\varphi(x)|\leqslant C\exp\{-\alpha^2x^2/4\}$ for all $x\in\bR$ (by \eqref{alp}, $\alpha^2<2$). By \eqref{1/alp} and \eqref{alp},
\[I_1=|\Phi(v/s^H)-\Phi(v/(\alpha s^H))|\leqslant \varphi\Big(\frac{v}{\alpha s^{H}}\Big)\frac{|v|}{\alpha s^{H}}\, \alpha(1-1/\alpha)\leqslant\]
\[\leqslant \alpha C \exp\{-v^2/(4s^{2H})\}\frac{(t-s)^H}{t^H}\leqslant  \sqrt{2} C \exp\{-v^2/(4t^{2H})\}\frac{(t-s)^H}{t^H}.\]

Let us now bound $\p(\alpha B_s<v<B_t)$. Recall that, for all $p>0$ and $q\in\bR$,
\[\int_{-\infty}^\infty \exp\{-px^2-qx\}dx=\sqrt{\frac{\pi}{p}}\,\exp\{q^2/(4p)\},\] and, as a result, 
 if $\xi\sim\mathcal N(0,1),$ then
\begin{equation}\label{alal}
\e\exp\{-(\sigma \xi+a)^2/2\}=\frac{1}{\sqrt{\sigma^2+1}}\exp\Big\{-\frac{a^2}{2(\sigma^2+1)}\Big\},\quad a\in\bR,\;\;\sigma>0.
\end{equation}
Using the latter  and \eqref{Phi}, we have
\begin{align*}
\p(\alpha B_s<v< B_t)&=\p(\alpha B_s<v< Z+\alpha B_s)\\
&=\e I(v>\alpha B_s) \Big[1-\Phi\Big(\frac{v-\alpha B_s}{\sqrt{t^{2H}-\alpha^2s^{2H}}}\Big)\Big]\\
&\leqslant \e \exp\Big\{-\frac{(v-\alpha B_s)^2}{2(t^{2H}-\alpha^2s^{2H})}\Big\}\\
&\quad=\e\exp\{-(\sigma \xi+a)^2/2\},
\end{align*}
where $\xi\sim\mathcal N(0,1),$ $a=v/\sqrt{t^{2H}-\alpha^2s^{2H}}$ and $\sigma=\alpha s^{H}/\sqrt{t^{2H}-\alpha^2s^{2H}}.$
For such $a$ and $\sigma$, 
\[\sigma^2+1=\frac{t^{2H}}{t^{2H}-\alpha^2s^{2H}}\quad\text{and}\quad \frac{a^2}{\sigma^2+1}=\frac{v^2}{t^{2H}}.\]
By  \eqref{ts}, $1/\sqrt{\sigma^2+1}=\sqrt{t^{2H}-\alpha^2s^{2H}}/t^H\leqslant (t-s)^H/t^H.$ Finally, applying \eqref{alal}, we conclude that $\p(\alpha B_s<v<B_t)\leqslant\exp\{-v^2/(2t^{2H})\} (t-s)^Ht^{-H}.$ Combining all bounds together, we get the desired inequality. Q.e.d.

\section*{Appendix A}
\noindent{\bf Proof of Lemma \ref{l1}.}
For simplicity, let $[t_0,t_n]=[0,1]$. Continue $f$ from  $\{(s,t):0\leqslant s\leqslant t\leqslant 1\}$ to  $\{(s,t):-1\leqslant s\leqslant t\leqslant 2\}$ by $f(s,t)=f(\varphi(s),\varphi(t))$, $s\leqslant t$, where $\varphi(t)=tI(0\leqslant t\leqslant 1)+I(t>1)$.  Then, for  $c=3d$ ($d<1/6$),
\[\iint\limits_{[-1,0]\times[0,c]}\frac{f(s,t)}{(s-t)^2}\,ds\,dt=
\int_{0}^{c} f(0,t)\,dt\int_{-1}^{0}\frac{ds}{(s-t)^2}
\leqslant\int_0^{c}\frac{f(0,t)}{t}dt,\] 
\[\iint\limits_{[c,1]\times[1,2]}\frac{f(s,t)}{(s-t)^2}\,ds\,dt=
\int_{c}^1 f(s,1)\,ds\int_1^2\frac{dt}{(s-t)^2}
\leqslant\int_{c}^1\frac{f(s,1)}{1-s}dt,\]
and  $f(s,t)=0$ if $s\leqslant t<0$ or $1<s\leqslant t$. Set  $g(s,t)=f(s,t)I(t<s+3d).$ Then
\[\int_{-1}^2\int_s^2\frac{g(s,t)}{(s-t)^2}\,ds\,dt\leqslant I_g(0,1).\]

We finish the proof by showing that 
\begin{equation}\label{DIc1-2}
\sum_{i=1}^n I_f(t_{i-1},t_i)\leqslant 9\int_{-1}^2\int_s^2 \frac{g(s,t)}{(s-t)^2}\,ds\,dt.
\end{equation}
Let further $0\leqslant u<v\leqslant 1$. By the monotonicity of $f$,
\begin{align*}
f(u,v)&=\frac1{\ln(4/3)}\iint\limits_{A(u,v)}\frac{f(u,v)}{(s-t)^2}\,ds\,dt\leqslant
4\iint\limits_{A(u,v)}\frac{f(s,t)}{(s-t)^2}\,ds\,dt.\end{align*} Here  $A(u,v)=[u-\Delta,u)\times[v,v+\Delta)$, $\Delta=v-u$. 

Let us show that
$A(t_{i-1},t_i)$ are pairwise disjoint. Fix $i<j$. If   $t_j\geqslant t_i+\Delta_i$  with $\Delta_i=t_{i}-t_{i-1}$, then $A(t_{i-1},t_i)\cap A(t_{j-1},t_j)=\varnothing$. If  $t_j\in[t_i,t_i+\Delta_i)$,  then \[t_{j-1}-\Delta_j\geqslant t_i-\Delta_j=t_i-t_j+t_{j-1}> -\Delta_i+t_i =t_{i-1} \] and, as a result, $A(t_{i-1},t_i)\cap A(t_{j-1},t_j)=\varnothing$.
 Therefore,  denoting the union of $A(t_{i-1},t_i)$ by $A$, we get
$$\sum_{i=1}^nf(t_{i-1},t_i)\leqslant9\iint\limits_{
A}\frac{f(s,t)}{(s-t)^2}\,ds\,dt.$$ In addition, setting
 $d_j=(v-u)/2^j$, we find that 
\[\int_u^v\frac{f(u,t)}{t-u}\,dt=
\sum_{j=1}^{\infty}\int_{u+d_j}^{u+d_{j-1}}\frac{f(u,t)}{t-u}\,dt
\leqslant\sum_{j=1}^\infty\int_{u+d_j}^{u+d_{j-1}}\frac{f(u,t)}{d_j}\,dt
\leqslant\]\[\leqslant\sum_{j=1}^\infty\;
\iint\limits_{[u-d_j,u)\times[u+d_j,u+d_{j-1})}
\frac{f(s,t)}{d_j^{\,2}}\,dt\,ds\leqslant
9\iint\limits_{B(u,v)}\frac{f(s,t)}{(s-t)^2}\,ds\,dt,\] where
\[B(u,v)=\bigcup_{j=1}^\infty\,\big([u-d_j,u)\times[u+d_j,u+d_{j-1})\big) =
\bigcup_{j=1}^\infty A\bigl(u,u+(v-u)/2^j\bigr).\] The sets 
$B(t_{i-1},t_i)$ are pairwise disjoint, since, for all $i$, \[B(t_{i-1},t_i)\subseteq \{(s,t):t\in[t_{i-1},t_i)\}.\] Denoting their union by $B$, we get 
\[\sum_{i=1}^n\int_{t_{i-1}}^{t_i}\frac{f(t_{i-1},t)}{t-t_{i-1}}dt\leqslant
9\iint\limits_{B}\frac{f(s,t)}{(s-t)^2}\,ds\,dt.\]
Since $A,B\subseteq C=\{(s,t)\in[-1,2]^2: s\leqslant t\leqslant s+3d\}$,
\[\sum_{i=1}^nf(t_{i-1},t_i)+\sum_{i=1}^n\int_{t_{i-1}}^{t_i}\frac{f(t_{i-1},t)}{t-t_{i-1}}dt\leqslant
13\iint\limits_{C}\frac{f(s,t)}{(s-t)^2}\,ds\,dt.\] Bounding the other terms  in $I_f(t_{i-1},t_i)$ in a similar way, we obtain \eqref{DIc1-2}. Q.e.d.

\end{document}